\documentclass[a4paper,twoside,ddraft,amssymb]{amsart}
\usepackage{amsmath}
\theoremstyle{plain}
\newtheorem{definition}{Definition}[section]

\newtheorem{lemma}[definition]{Lemma}
\newtheorem{theorem}[definition]{Theorem}
\newtheorem{corollary}[definition]{Corollary}
\newtheorem{proposition}[definition]{Proposition}

\newtheorem{notation}[definition]{Notation}
\newtheorem{remark}[definition]{Remark}
\theoremstyle{definition}

\theoremstyle{remark}

\begin{document}
\def\ZZ{Z\!\!\! Z}
\def\RR{I\!\! R}
\def\CC{I\!\!\! C}
\def\eop{\hfill QED.\vskip 0.3cm }
\def\cht{\check\tau}
\def\chs{\check\sigma}
\def\chr{\check\rho}
\def\nn{{a_1\over | X |}}
\def\n{\omega }
\def\p #1{p_{1,1}^{#1}}
\def\L{{\mathcal{L}}}
\def\La{{\mathbf{L}}}
\def\k{{d}}
\def\gcud{\cdot \!\!>}
\def\lcud{<\!\!\cdot}
\def\gcunor{\cdot \!\!\!>}
\def\lcunor{<\!\!\!\cdot}
\def\gcu{\ \cdot \!\!\!\!>}
\def\lcu{<\!\!\!\!\cdot\ }
\def\rk{{\rm rk}}
\def\qbin #1#2{{#1\brack #2}_{\!\! q}}
\def\bperp{{\perp\!\!\!\!\perp }}
\def\FA{\ \forall}
\def\direct{\mathrel\triangleright\joinrel\mathrel<}
\def\vaca{(q^{d-1}-1)}
\def\eldelta{{1_{_X}}}
\def\direct{\mathrel\triangleright\joinrel\mathrel<}

\title{Tight frames for  eigenspaces of the Laplacian on dual polar graphs}
\author{F. Levstein*, C. Maldonado* and D. Penazzi*}
\address{* Universidad Nacional de C\'ordoba,
Facultad de Matem\'atica, Astronom\'\i a y F\'\i sica, C\'ordoba,
Argentina, Haya de la Torre y Medina Allende, +54-351-4334051/363.
CIEM-CONICET}
\thanks{Partially supported by Secyt-UNC, CONICET, ANPCyT}
\maketitle

\begin{abstract}
We consider $\Gamma=(X,E)$ a dual polar graph and we give a tight
frame on each eigenspace of the Laplacian operator associated to
$\Gamma$. We compute the  constants associated to each tight frame
and as an application we give a formula for the product in the
Norton algebra attached to the eigenspace corresponding to the
second largest eigenvalue of the Laplacian.
\end{abstract}

\section{Introduction}\label{dpg}\noindent

In algebraic combinatorics a lot of research has been done on
distance regular graphs. The main examples are the following
families: Johnson, Grassmann, Hamming and dual polar graphs.

In this paper we consider  the set of functions $\RR^{X}=\{ f:X
\rightarrow \RR\}$ where $X$ is the set of vertices of the dual
polar graphs. The distance on the graph gives rise to a Laplacian
operator $\L$ on $\RR^{X}$ and a decomposition of $\RR^{X}$ into
eigenspaces of $\L$. These topics can be seen in \cite{DR,S1,S2}.

First  we associate a lattice to the graph and characterize the
eigenspaces of $\L$ in terms of this lattice. Instead of an
orthogonal basis we can give a tight frame for each eigenspace.
The theory of finite normalized tight frames has seen many
developments and applications in recent years. See for instance
the references in \cite{BF, FMRHO, KC1, KC2, MFRHO, VW1, VW2}.

The eigenspace corresponding to the second largest eigenvalue of
$\L$ is of particular importance since one can reconstruct the
whole graph from the projections of the canonical basis onto it.
We explicitly compute the constant of the tight frame attached to
this eigenspace.

The notion of Norton algebra was developed to give realizations of
the finite simple groups as automorphisms group of an algebra. The
general construction starts with a graded algebra
$\mathcal{V}=\bigoplus_i \mathcal{V}_i$  and gives an algebra
structure on each subspace $\mathcal{V}_i$ by multiplying on
$\mathcal{V}$ and then projecting onto $\mathcal{V}_i$.

As an application we answer a problem posed to us by Paul
Terwilliger: a formula for the product in the Norton algebra
attached to the eigenspace corresponding to the second largest
eigenvalue of $\L$.

This article is organized as follows: In section \ref{def} we give
some classical definitions. In section \ref{ldpg},  we associate a
lattice to a dual polar graph $\Gamma$. In section
\ref{embedding}, we give  a convenient description for the
eigenspaces $V_i$ of $\L$. In the next section,  Theorem
\ref{tightframe} gives a tight frame on each eigenspace $V_i$ and
give a formula for the constant associated.

In the last section we compute an explicit formula for the product
in the Norton algebra mentioned above.

\section{Definitions}\label{def}
\subsection{Distance regular graphs}\label{drg}\noindent

Given  $\Gamma=(X,E)$ a graph with distance  $d(\ ,\ )$  we say
that it is distance regular if for any $(x,y) \in X\times X$ such
that $d (x,y)=h$ and for all $i,j\ge 0$ the cardinal of the set
$$\{z \in X \mid
d (x,z)=i \ \mbox{and}\ d (y,z)=j \}$$  is a constant denoted by
$p_{ij}^h$ which is independent of the pair $(x,y).$

\subsection{Adjacency algebra of a distance regular graph}\noindent

Let  $\Gamma=(X,E)$ be a distance regular graph of diameter $d$.
Let $Mat_{X}(\RR)$ denote the $\RR$-algebra of matrices with real
entries, where the rows and columns are indexed by the elements of
$X$.\\ For $0 \leq i \leq d$, let $A_i$ denote the following
matrix in $Mat_{X}(\RR)$:
\[ (A_i)_{xy} = \left \{ \begin{array}{ll}
1 & \mbox{if} \ d(x,y)=i \\
0 & \mbox{if} \ d(x,y)\neq i
\end{array}
\right. \] We call $A_i$ the {\it ith adjacency matrix} of
$\Gamma$. Using the definition it is not difficult to prove that
the adjacency matrices of a distance regular graph satisfy:

 (i') $A_0=I$ where $I$ is the identity matrix in
$Mat_{X}(\RR)$;

(ii') $A_0+\dots +A_d=J$ where $J$ is the all $1' s$ matrix in
$Mat_{X}(\RR)$;

(iii') $ A_iA_j=\sum_{h=0}^d p_{ij}^h A_h \ (0\leq i,j \leq d)$.;

(iv') ${A_i}^t=A_i$ \\
It follows from (i')-(iv') that $A_0,\dots, A_d$ form a basis for
a subalgebra $\mathcal{A}$ of $Mat_{X}(\RR)$. We call
$\mathcal{A}$ the {\it adjacency algebra} of $\Gamma$.

\vspace{2em}

It is known that the space of functions $\RR^{X}=\{ f:X
\rightarrow \RR\}$ has a decomposition $$\RR^{X}=\oplus_{j=0}^{d}
W_j$$ where $\{W_j\}_{j=0}^{d}$ are the common eigenspaces of
$\{A_i\}_{i=0}^{d}$. Let $p_i(j)$ the eigenvalue of $A_i$ on the
eigenspace $W_j$.

By Proposition 1.1 of section 3.1 of Chapter III of \cite{BI},
 the adjacency matrices $\{A_i\}_{i=0}^{d}$ and the
eigenvalues $\{p_i(j)\}_{i,j=0}^{d}$ of a given distance regular
graph $\Gamma$  satisfy: $$ A_i=v_i(A_1) , \quad
p_i(j)=v_i(\theta_j)$$ where $\theta_j=p_1(j)$, and
$\{v_i\}_{i=0}^{d}$ are polynomials of degree $ i$.

 We will order
the decomposition according to $\theta_0 > \theta_1
>...>\theta_d$.

In Theorem 5.1 of III.5 of \cite{BI}, one can find formulas for
the polynomials associated to each $\Gamma$.

\subsection{Dual Polar Graphs}\label{dpg}\noindent

Let $V$ be a finite dimensional vector space over $GF(q)$ (the
finite field with $q$ elements), together with a nondegenerate
form $\omega$. A subspace of $V$ is called {\it isotropic}
whenever the form vanishes completely on it. The dual polar graph
corresponding to $(V,\omega)$ is the graph $\Gamma=(X,E)$ where
\begin{eqnarray*}
X &=& \{ v \subseteq V : \ v\ \mbox{is maximal isotropic subspace} \} \\
E &=&\{(u,v) \in X\times X \ : \  ,\ dim(u\cap v))=d-1\}
.\end{eqnarray*}

\vspace{1em}

The dual polar graphs are distance regular, and are listed in page
274 of  \cite{BCN}. They are the following:
\begin{eqnarray*}
C_d(q): V&=&GF(q)^{2d}\ \mbox{with a nondegenerate sympletic form.}\\
B_d(q): V&=&GF(q)^{2d+1}\ \mbox{with a nondegenerate quadratic form.}\\
D_d(q): V&=&GF(q)^{2d}\ \mbox{with a nondegenerate quadratic form
of Witt index}\ d.\\
 ^2D_{\! d+\!1}(q): V&=&GF(q)^{2d+\!2}\ \mbox{with a
nondegenerate quadratic form of Witt index}\ d.\\
^2A_{2d}(r): V&=& GF(r^2)^{2d+1}\ \mbox{with a nondegenerate
Hermitean
form.}\\
^2A_{2d-1}(r): V&=&GF(r^2)^{2d}\ \mbox{with a nondegenerate
Hermitean form.}
\end{eqnarray*}

\vspace{.2cm}

In each of the cases above, the dimension of the maximal isotropic
spaces is ${d}$.

We will denote $U^{\bperp }=\{v\in V:\omega(v,u)=0\FA\ u\in U\}$.

In each case there is a group acting on these spaces, namely the
group $G_\omega$ of  linear transformations on the underlying
space $V$ that preserve the form $\omega$.

\section{Lattice associated with dual polar graphs}\label{ldpg}

In this section we consider the graphs defined above and we
associate a lattice to them. We recall the following definitions:

\begin{itemize}
\item A partial order is a binary relation "$\leq$" over a set P
which is reflexive, antisymmetric, and transitive.
\item A partially ordered set (POSET) $(P,\leq)$ is a set $P$ with a partial order $\leq$.
\item A lattice $(\La,\leq,\wedge,\vee)$ is a POSET  $(\La,\leq)$ in
which  every pair of elements $u, w \in \La$ has a least upper
bound and a greatest lower bound. The first is called  the join
and it is denoted by $u\vee w$ and the second is  called the meet
and it is denoted by $u\wedge w$.

\end{itemize}

\subsection{Construction of the lattice.}\noindent

Let $\Gamma=(X,E)$ be a dual polar graph and be $V$ be the
underlying finite dimensional vector space over $GF(q)$.
\begin{definition}\noindent
\begin{eqnarray*}
& \Omega_\ell&=\{v \subseteq V : \ v\  \mbox{is an isotropic
subspace and}\
dim(v)=\ell \}.\ \ell=0,... d. \\
 & \Omega_{d +1}&=\{V\} \end{eqnarray*}

We let $\ \hat{0}:=\{0\}$ and $\ \hat{1}:=V$  and  we denote
$\Omega_d$ by $X$.

\end{definition}

 We will always work with $d>1$, i.e., $\Omega_1 \ne X$.

\begin{definition}\noindent\label{la}

 \begin{itemize}
 \item $\La=\cup_{\ell=0}^{\k +1}\Omega_\ell.$
\item Given isotropic subspaces $u, w \subseteq V $ we set: \begin{itemize}
 \item  $u\le w$ if and only if, $u$ is a subspace of $w$.
 \item  $u\wedge w=u\cap w$
 \item  $u\vee w=span\{u,w\}$ if that space is isotropic, otherwise, $u\vee w=V=\hat{1}$ \end{itemize}
\item The rank of $w\in \Omega_\ell$ is $\ell$ and it is denoted by $\rk (w)$.
\item Given $w\in \Omega_j$; $u$ covers $w$ or $w$ is
covered by $u$, if $u\in \Omega_{j+1}$ and $w\le u$. We denote it
by $u\gcunor w$ or $w\lcunor u$.
\item An atom is an element that covers $\hat 0$ and a coatom is an element
covered by $\hat 1$.

 \end{itemize}

 \end{definition}

Is not difficult to see that  $(\La, \leq, \wedge , \vee)$ is a
finite lattice with lowest element $\hat 0$ and greatest element
$\hat 1$. In our notation, the set of atoms is $\Omega_1$ and the
set of coatoms is $X$.

\vspace{.5em}

\begin{lemma}\label{prop}

The lattice $\La$  has the following properties:
\begin{enumerate}
    \item $\La$ is atomic.
    \item $ u\vee w\not=\hat 1\Rightarrow \rk(u)+\rk(w)=\rk(u\vee w)+\rk(u\wedge
    w)$\label{modularity}
\end{enumerate}
\end{lemma}

\proof \noindent

\begin{enumerate}

\item Each element $u \in \Omega_j$ of the lattice is a subspace of $GF(q)^n$, so
taking a basis $\{v_1,...,v_j\}$ of $u$, we obtain that
$u=span(v_1)\vee span(v_2)\vee ...\vee span(v_j)$ is a join of
atoms.
\item The rank of an element is the dimension, so the formula is
true because of the well known identity
$dim(u+w)=dim(u)+dim(w)-dim(u\cap w)$. (The formula fails for the
case $u\vee w=\hat1$ because then $u\vee w$ is not equal to
$u+w$).\end{enumerate} \eop
\begin{corollary}\label{first}
If $\tau $ and $\sigma $ are different atoms such that \
$\tau\vee\sigma\ne \hat 1 $, then \\$$\rk(\tau\vee\sigma)=2 .$$
\end{corollary}
\proof
\begin{eqnarray*}
\rk(\tau \vee \sigma)&=&\rk(\tau)+\rk(\sigma)-\rk(\tau\wedge\sigma)\\
&=&1+1-0\\
&=&2\\
\end{eqnarray*}
\eop

\begin{lemma}\label{coverproperty}
Let $u$ and $w$ be elements of the lattice which are not coatoms. If
 $u\vee w$ covers both $u$ and $w$ then $u$ and $w$ both cover $u\wedge
 w$.

 Reciprocally if $u$ and $w$ cover $u\wedge
 w$ and $u\vee w\ne\hat 1$, then $u\vee
 w$ covers  both $u$ and $w$.

\end{lemma}
\proof In order to prove the first statement, observe that  $z$
covers $w$ iff $z\ge w$ and $\rk(z)=\rk(w)+1$. So, $u\vee w$
covers both $u$ and $w$ iff $\rk(u\vee w)=\rk(u)+1=\rk(w)+1$ (in
particular, we must have that $\rk(u)=\rk(w)$). Also, since $u$
and $w$ are not coatoms and $\rk(u\vee w)=\rk(u)+1$ we deduce that
$u\vee w\ne\hat 1$. Then, by  Lemma \ref{prop}
(\ref{modularity}),we get $\rk(u)+\rk(w)-\rk(u\wedge w)=\rk(u)+1$,
i.e., $\rk(w)=\rk(u\wedge w)+1$, which implies that $w$ covers
$u\wedge w$. The proof is similar for $u$.

 Reciprocally, if $u$ and $w$ cover $u\wedge
 w$, then $\rk(u)=\rk(w)=\rk(u\wedge w)+1$. Using  Lemma \ref{prop}
(\ref{modularity}) we get $\rk(w)=\rk(u)+\rk(w)-\rk(u\vee w)+1$
which implies $\rk(u)+1=\rk(u\vee w)$ and then that $u\vee w$
covers $u$ (and similarly $w$). \eop

\section{Description of the eigenspaces of a dual polar graph using the associated lattice}\label{embedding}

\vspace{1em}

In this section we will consider a dual polar graph
$\Gamma=(X,E)$, together with its associated decomposition:
$$\RR^{X}=\oplus_{i=0}^{d} W_i ,$$ where $\{W_i\}_{i=0}^{d}$ are the
common eigenspaces of the adjacency matrices of $\Gamma$.

We will describe each of the eigenspaces $\{W_i\}_{i=0}^{d}$,
using the lattice previously defined.

\vspace{1em}

For ease of writing, we will use the following notation:
\begin{notation}\noindent
\begin{itemize}
    \item For any statement $P$, let $[P]=\begin{cases}  1 &
\mbox{if}\  P\ \mbox{is true.}\ \cr 0& \mbox{if } P \ \mbox{is
false.}\
\end{cases}
$
    \item $\RR^X=\{f:X \rightarrow \RR \}$,\\ (thus $0\in  \RR^X$
    will denote the function such that $0(x)\equiv 0 , \forall x \in
    X$,\\ analogously $1 \in \RR^X$.)
    \item Let $<;>$ be the inner product in $\RR^X$ defined by $<f,g>=\sum_{x\in
X}f(x)g(x)$.
\item For $u\subseteq  \RR^X$ let $u^{\perp}=\{f \in \RR^X :\ <f;g>=0 \ \forall \ g \ \in U\}$
    \item $||f||^{2}=<f;f>$.
\end{itemize}

\end{notation}

We will need the following lemma.

Recall that $\qbin {i} 1=\left\{
\begin{array}{ccc} \frac{q^{i}-1}{q-1} & \forall \
i \geq 1\\ 0 & \forall \ i < 1 \end{array}\right .$ and $\qbin {i}
j=\frac{\qbin {i} 1 \qbin {i-1} 1 ... \qbin {i-j+1} 1}{\qbin {j} 1
\qbin {j-1} 1 ... \qbin {1} 1}$.

\begin{lemma} \label{BCN} (9.4.2 of \cite{BCN})\noindent

Let $e$ be $1,1,0,2,{3\over 2},{1\over 2}$ in the respective
cases:

$C_d(q),B_d(q), D_d(q), ^2D_{d+1}(q), ^2A_{2d}(r), ^2A_{2d-1}(r)$.

Let $W$ be a fixed isotropic space of dimension $j$. The number of
isotropic spaces $U$ of dimension $(k+l+m)$, meeting $W$ in a
space of dimension $m$ and $W^{\bperp } $ in space of dimension
$l+m$ is:
\begin{eqnarray*}
q^{l(j-m)+k(2d-j-m-2l+e-1)-k(k-1)/2}\qbin jm\qbin {j-m}k\qbin
{d-j}l\prod_{i=0}^{l-1}(1+q^{d+e-j-i-1})
\end{eqnarray*}

\end{lemma}

\begin{lemma}\label{al}

For $z\in\Omega_j$, let $a_j=\vert \{x\in X:z\le x\} \vert$. Then
$$a_j=
\prod_{i=0}^{d-j-1} (1+q^{e+i}) \quad {\rm if\ } 0\leq j  \leq d
 \quad   {\rm and } \quad  a_{d+1}=0$$
\end{lemma}
\proof

We use the previous lemma with $W=z \in \Omega_j$ and $U=x \in
\Omega_d$. Then $k+l+m=d$, $m=j$ and $l+m=d$. Therefore we can
apply the formula  with $k=0,\ l=d-j$. \eop

\vspace{.5em}

\begin{remark}\label{remal} Note that $a_j=(1+q^{d-j-1+e})\ a_{j+1}$.
\end{remark}

\vspace{.5em}

\begin{definition} \noindent

$\iota:\La  \rightarrow  \RR^X$ is the map defined by $
\iota(z)(x)= [z\le x] \ \forall \ z \in \La ,\ x\in X $

\end{definition}

\begin{lemma}\label{iotahat}\noindent

$$i)\ \iota(\hat 1)=0 \ \in \ \RR^{X} \qquad ii)\ \iota(\hat 0)=1 \ \in \ \RR^{X} \qquad iii)\
\iota(z)\iota(y)=\iota(z\vee y) \ \forall \ z,y \in \La .$$

\end{lemma}
\proof \noindent

$i) \iota(\hat 1)(x)=[\hat 1\le x]=0\FA x$ , since $\hat 1$ is
above the $x$'s.

$ii) \iota(\hat 0)(x)=[\hat 0\le x]=1\FA x$.

$iii)$
 \begin{eqnarray*}
\iota(z)(x)\iota (y)(x)&=&[z\le x][y\le x]\\
&=&[(z\le x) \ and \  (y\le x)]\\
&=&[z\vee y\le x]\\
&=&\iota(z\vee y)(x)
\end{eqnarray*}
\eop
\begin{lemma}\label{innerproduct}
$$<\iota(z);\iota(y)>=||\iota(z\vee y)||^2 \quad \forall \ z,y \in \La .$$

\end{lemma}
\proof
\begin{eqnarray*}
<\iota(z);\iota (y)>&=&\sum_{x\in X}\iota(z)(x)\iota (y)(x)\\
&=&\sum_{x\in X}\iota(z\vee y)(x)\quad ({\rm by\ Lemma\ \ref{iotahat}\ iii)})\\
&=&\sum_{x\in X}(\iota(z\vee y)(x))^2\quad ({\rm since\ }\iota(z\vee y)(x)\in\{0,1\} ) \\
&=&||\iota(z\vee  y)||^2
\end{eqnarray*}
\eop
\begin{corollary} $z\vee y=\hat 1$ if and only if $\iota(z)$ and $\iota(y)$ are orthogonal to each other.
\end{corollary}
\proof By lemmas \ref{innerproduct} and \ref{iotahat} i),
$<\iota(z);\iota(y)>=0$ if and only if $z\vee y=\hat 1$. \eop

\begin{lemma}\label{norm}
$$||\iota(z)||^2=a_j\ \FA \ z\in \Omega_j, \ j=0,1,...,d$$
\end{lemma}
\proof
\begin{eqnarray*}
||\iota(z)||^2&=&\vert\{x\in X : \iota(z)(x)=1\}\vert\\
&=&\vert\{x\in X : z \leq x\}\vert\\
&=&a_j
\end{eqnarray*}
\eop
\begin{corollary}\label{innerproductcoro}
If $z\vee y\in \Omega_j$, then $<\iota(z);\iota (y)>=a_j$
\end{corollary}
\proof Direct from the two previous lemmas.\eop

\begin{lemma}\label{iotainner}
If $\tau $ and $\sigma$ are both atoms then:
$$<\iota(\tau );\iota(\sigma)>=\begin{cases}
a_1& {\rm if}\ \tau=\sigma \cr 0& {\rm if}\ \tau \vee
\sigma=\hat{1} \cr a_2& {\rm otherwise}
\end{cases}$$
\end{lemma}
\proof

If $\tau=\sigma$, then $\tau \vee \sigma=\tau $, and so we have
that $<\iota(\tau );\iota (\sigma)>=a_1$ by Lemma \ref{norm}. If
$\tau\vee\sigma=\hat 1$, then by Lemma \ref{innerproduct}, we have
that $<\iota(\tau );\iota(\sigma)>=0$.

If $\tau\ne \sigma$ and $\tau\vee\sigma\ne\hat 1$, then by Lemma
\ref{first}, $\tau\vee\sigma\in \Omega_2$, so by Corollary
\ref{innerproductcoro}, we have $<\iota(\tau
);\iota(\sigma)>=a_2$. \eop

\subsection{A filtration for $\RR^{X}$.}

\begin{definition}\label{deflambdaj}\noindent

For $j=0,1,...,d$, let $\Lambda_j \subseteq \RR^{X}$ be the
subspace generated by $\{\iota(x)\}_{x \in \Omega_j}$, that is
$\Lambda_j=span(\iota(\Omega_j))$.
\end{definition}
We want to show that $\Lambda_j\subseteq\Lambda_{j+1}$. For this
we need some tools first.
\begin{definition}
 Given $w\in \La$, let:
 $$w^*=\sum_{v\gcud w}\iota(v)$$
\end{definition}
\begin{lemma}\label{lemawstar}
Given $w\in \La$, $\iota(w)$ is a scalar multiple of $w^*$. In
fact,
$$w^*=\qbin {{d}-j}1\iota(w) \quad \forall\ w \in \Omega_j$$
\end{lemma}
\proof
Given $x\in X$, we have:
 \begin{eqnarray*}
w^*(x)&=&\sum_{v\gcud w}[v\le x]\\
&=&\vert\{v:w\lcunor v\le x\}\vert \quad (*)\\
\end{eqnarray*}
Clearly, if $w\not\le x$, that number is zero, i.e.,
$w^*(x)=0=\iota(w)(x)$ if $w\not\le x$. On the other hand, if
$w\le x$, the number in $(*)$ is the number of spaces in $x$ built
from $w$ by adding a one-dimensional space. That one-dimensional
space must be in $x$ and not in $w$, so there are $\qbin {{d}-j}1$
ways of doing this. So
$$ w^*(x)=\left\{\begin{array}{ccc r}
0  & {\rm if}& w\not\le x &\\
\qbin {{d}-j}1 & {\rm if}&  w\le x ,& {\rm then}
\end{array} \right.$$
$$w^*=\qbin {{d}-j}1 \iota(w).$$\eop

\begin{corollary}\label{corolambdas}
$$\Lambda_0\subseteq\Lambda_1\subseteq...\subseteq \Lambda_{d}=\RR ^X$$
\end{corollary}
\proof Let $f\in \Lambda_j$. We can assume that $f=\iota(w)$. By
definition, $w^*\in \Lambda_{j+1}$. But by lemma \ref{lemawstar},
$\iota(w)$ is a non-zero  scalar multiple of $w^*$,
 so $f\in\Lambda_{j+1}$

\eop
\begin{definition}
Let   $V_0=\Lambda_0$ and $V_j=\Lambda_j\cap\Lambda_{j-1}^\perp
\quad j=1,...,d$.
\end{definition}
So, we have that $\Lambda_j=V_0\oplus V_1\oplus ...\oplus V_j$.

We want to show that for $j=1,...,d$, $V_j\neq  \{0\}$, that is
$\Lambda_{j-1} \neq \Lambda_{j}$ . To prove this, we need more
lemmas.
\begin{definition}

Let $\L:\RR^X\mapsto \RR^X$ denote  the  Laplacian operator
defined by $$\L(f)(x)=\sum_{y\in X :\ d(x,y)=1}f(y)$$
\end{definition}

Observe that
\begin{eqnarray*}
\L (f)(x)&=&\sum_{y\in X}[d(x,y)=1]f(y)\\&=&\sum_{y\in X}(
A_1)_{xy}f(y) ,\end{eqnarray*} where $A_1$ is the first adjacency
matrix of $\Gamma=(X,E)$ a dual polar graph. So $\L$ can be
thought as multiplication by $ A_1$. In particular, $\L$ is
symmetric and $<\L(f),g>=<f,\L(g)>$.
\begin{lemma}
If $x\in X$, then $\L(\iota(x))=\sum_{y\in X :\
d(x,y)=1}\iota(y)$.
\end{lemma}
\proof
\begin{eqnarray*}
\L(\iota(x))(z)&=&\sum_{y\in X:\ d(z,y)=1}\iota(x)(y)\\
&=&\sum_{y\in X}[d(z,y)=1][x\le y]\\
&=&\sum_{y\in X}[d(z,y)=1][x=y]\qquad {\rm because\ } x\le y\iff x=y\ {\rm since\ } x,y\in X\\
&=&[d(z,x)=1]\ \ {\rm while} \\
\\
(\sum_{y\in X:\ d(x,y)=1}\iota(y))(z)&=&\sum_{y\in X}[d(x,y)=1][y\le z]\\
&=&\sum_{y\in X}[d(x,y)=1][y=z]\\
&=&[d(x,z)=1]\\
\end{eqnarray*}
\eop
\begin{lemma}\label{lemasumaiotas}
Let $x\in X$. Then:
$$\L(\iota(x))=-\qbin {d} 1\iota(x)+\sum_{z: z\lcud x}\iota(z)$$
\end{lemma}
\proof

\begin{eqnarray*}
(\sum_{z: z \lcud x}\iota(z))(y)
&=&\sum_{z: z \lcud x}[z\le y]\\
&=&\sum_{z\in \Omega_{{d}-1}}[z\le x\wedge y]\\
&=&\begin{cases}
 \sum_{z\in \Omega_{{d}-1}}[z\le x]& if x=y \cr
\sum_{z\in \Omega_{{d}-1}}[z=x\wedge y]& if x\ne y \end{cases}\\
&=&\begin{cases} \qbin {d} 1& if x=y\cr [x\wedge y\in
\Omega_{{d}-1}] &if x\ne y \end{cases}\\
&=&\qbin {d} 1[x=y]+[x\wedge y\in \Omega_{{d}-1}] [x\ne y]\\
&=&\qbin {d} 1[x=y]+[d(x,y)=1]\\
&=&\qbin {d} 1\iota(x)(y)+\L(\iota(x))(y)\\
\end{eqnarray*}
\eop

\def\DDelta{(q^j+q^{d+e-j-1})}
\begin{proposition}\label{arribayabajo}\noindent

Let $\Gamma=(X,E)$ be a dual polar graph  and let $e$ be as in
Lemma \ref{BCN}.

For  $j=0,1,...,d-1$ and  for all $w\in \Omega_j$ the following
holds:
$$\sum_{u:u\gcud w}\sum_{z: z\lcud u}\iota(z)=\qbin {d-j}1\left(
\DDelta\iota(w)+\sum_{v: v\lcud w}\iota(v)\right)
$$
\end{proposition}
\proof Let us call $S(x)$ the function on the left hand side of
the equation above and $R(x)$ the function on the right  hand
side.

We have to see that evaluating on an arbitrary $x\in X$, they are
both equal.
\begin{enumerate}

\item
{\bf CASE 1}: $w \wedge x \in \Omega_j$.

\vspace{.5em}

In this case, $w \wedge x =w$ that is  $w\leq x$. Then
$$R(x)=\qbin{d-j}1\left(\DDelta+ \vert\{v :\ v\lcu w \ {\rm and}\
v\le x\}\vert\right)$$
 However, $v\lcu w\le x  \Rightarrow v\le x$, so
$$R(x)=\qbin{d-j}1\left( \DDelta+ \vert\{v :\ v\lcu w\}\vert \right)$$

 Since $\vert\{v :\ v\lcu w\}\vert$ is the number of spaces of dimension
$j-1$ in a  space of dimension $j$, i.e., $\qbin j{j-1}=\qbin j1$, we conclude:

\begin{eqnarray*}
R(x)&=& \qbin{d-j}1\left(\DDelta+\qbin j1\right)
\end{eqnarray*}

On the other hand, $S(x)=\sum_{u\gcud w}\vert\{z :\ z\lcu u \ {\rm
and} \ z\le x\}\vert$(recall that $w, z \in \Omega_j , \ u \in
\Omega_{j+1}$ and $x \in \Omega_d$).

\vspace{.5em}

If $u\le\! x$ then $z\lcu u\!\Rightarrow\! z\le\! x$, i.e.:
 $$\vert\{z: \ z\lcu u \ {\rm
and} \ z\le x\}\vert=\vert\{z: \ z\lcu u\}\vert=\qbin {j+1}j=\qbin
{j+1}1$$

If $u\not\le x$, then $x\wedge u\ne u$. But $w\le x\wedge u\le u$
and $w\lcu u$, so $x\wedge u=w$.

So in this case $\{z: \ z\lcu u : z\le x\}= \{w\}$ and thus,

$\vert\{z: \ z\lcu u : z\le x\}\vert=1$. Therefore,
\begin{eqnarray*}
S(x)&=&\vert\{u : \ w\lcu u\le x\}\vert\qbin {j+1}1+\vert\{u: \
u\gcu w \ {\rm
and} \  u\not\le x\}\vert\\
&=&\qbin{d-j}1\qbin {j+1}1+\vert\{u: \ u\gcu w \ {\rm and} \
u\not\le x\}\vert
\end{eqnarray*}

This  last number can be computed from $$|\{u: u\gcu w \ {\rm and}
\ u\not\le x\}|= |\{u :\ u\gcu w\}| - |\{u :\  x\geq u\gcu w\}|$$

For this, we need Lemma \ref{BCN}.

In order to compute $|\{u :\ u\gcu w\}| $ we fix $w$ isotropic of
dimension $j$ and we want to compute the number of isotropic
spaces $u$ of dimension $j+1$. Since $u \geq w$ we have $dim (u
\cap w)=dim\ w=j$ and since $u$ must be isotropic, $u\subseteq
w^\bperp$ and then $dim (u\cap w^\bperp)=dim\ u=j+1$. Then in the
notation of the lemma ``$W$"$=w$, ``$U$"$=u$ and
$$ k+l+m=j+1, \quad m=j, \quad
{\rm and}\quad l+m=j+1$$ Therefore $l=1$, $k=0$.

Then, the lemma gives us
\begin{eqnarray*}
|\{u : u\gcu w\}|&=& \qbin {d-j}1 (1+q^{d+e-j-1})\\
|\{u: u\gcu w \ , \ u\not\le x\}|&=& |\{u\gcu w\}| - |\{u : x \geq u\gcu w\}|\\
&=&\qbin {d-j}1 (1+q^{d+e-j-1}) - \qbin
{d-j}1 = \qbin {d-j}1 q^{d+e-j-1} \\
{\rm and} \ S(x)&=&\qbin{d-j}1\left(\qbin
{j+1}1+q^{d+e-j-1}\right)\\ \Rightarrow  S(x)&=&R(x)
\end{eqnarray*}

\medskip

\item
{\bf CASE 2}:  $w\wedge x\in \Omega_k,$ for  $k< j-1$.

\vspace{.5em}

In this case, $w \not\le x$ and  $w\wedge x\not\lcu w$. Then
$S(x)$ is the number of pairs $(u,z)$ which satisfy:
\begin{enumerate}
\item $w\lcu u$
\item $z\lcu u$
\item $z\le x$
\end{enumerate}
(Recall that $w, z \in \Omega_j , \ u \in \Omega_{j+1}$ and $x \in
\Omega_d$).

Now, $((c)+w\not\le x)\Rightarrow z\ne w$. This plus $(a)$ and
$(b)$ implies $w\vee z=u$, i.e., $w,z\lcu
w\vee z$. Then, $w\wedge z\lcu w,z$ (Lemma \ref{coverproperty}).\\
However, $(c)\Rightarrow w\wedge z\le w\wedge x \ (\le w)$, so
since $w\wedge z\lcu w$, we have only two options for $w\wedge
x$:\begin{itemize}
    \item $w\wedge x=w\wedge z$ (impossible since $w\wedge x\not\lcu w$ while $w\wedge z\lcu w$)
    \item  $w\wedge x=w$, which would imply $w\le x$.
\end{itemize}     Since neither of
these happens in CASE 2, we conclude that there are no such pairs
$(u,z)$, and $S(x)=0$.

On the other hand, $$R(x)=\qbin{d-j}1\left( \DDelta\ 0+ \vert\{v:\
v\lcu w\ {\rm and } \ v\le x\}\vert \right )$$ But $(v\lcu w\ {\rm
and } \ v\le x)\Rightarrow v\le w\wedge x\le w$.

Again, $w\wedge x\neq w$, because $w\not\le x$ and  $w\wedge x\neq
v$ because $v\lcu w$ while $w\wedge x\not\lcu w$. So $R(x)=0$ too
and we have equality in this case.

\medskip

\item
{\bf CASE 3}: $w\wedge x\in \Omega_{j-1}$

In this case $w \not\le x$ and  $w\wedge x\lcu w$.

As in CASE 2, $R(x)=\qbin {d-j}1 \vert\{v:\ v\lcu w\ {\rm and } \
v\le x\}\vert$, and  given such a $v$, we deduce that either
$w\wedge x=w$ (impossible since $w \wedge x\in \Omega_{j-1}$) or
$w\wedge x=v$.  Then $\{v: \ v\lcu w\ {\rm and } \ v\le
x\}=\{w\wedge x\}$  so $
 R(x)=\qbin {d-j}1 .$\\
On the other hand $$S(x)=\sum_{u: \ u\gcud w}\sum_{z}\ [z\lcu
u][z\le x]=\sum_{u: \  u\gcud w}\sum_{z}\ [z\le u\wedge x][z\lcu
u]$$ But $z\le u\wedge x\le u$ and $z\lcu u$ imply that $u\wedge
x$ is either $z$ or $u$.\\ If $u\wedge x=u,\  \Rightarrow u\le x$,
and since $w\lcu u$, then $w\le x$, which can not happen in CASE
3. Therefore, $u\wedge x=z$, i.e.,
$$\sum_{z}\ [z\le u\wedge x][z\lcu
u]=\left\{\begin{array}{ccc}  1 & \ \mbox{if}\ u\wedge x\lcu u \\
 0 & \ \mbox{otherwise}\end{array}\right.
$$

So assume  $ u\wedge x\lcu u$. Since $u\wedge x\le x$ and
$w\not\le x$, we must have $u\wedge x\ne w$. Since $u\gcu w$ and
$u\gcu u\wedge x$, then $u=(u\wedge x)\vee w$, i.e., $u$ is
determined by $z=u\wedge x$, i.e.,
$$S(x)=\vert\{z:\ z\le x\ {\rm and } \ z\vee w\gcu w,z\}\vert.$$ Since $z\vee w\gcu
w,z\Rightarrow w\wedge z\lcu w$ and $w\wedge z\le w\wedge x\lcu
w$, we conclude that $w\wedge z=w\wedge x$. Also, in order to have
$z\vee w\gcu w$, i.e., to be isotropic, we must have $z\le
w^\bperp$. Hence, $S(x)=\vert\{z:\ w\wedge x\lcu z\le w^\bperp\cap
x\}\vert=\qbin {\theta(x)}1$, where $\theta(x)=dim(w^\bperp\cap
x)-dim(w\wedge x)=dim(w^\bperp\cap x)-(j-1)$.

Since $w\wedge x\lcu w$, we have that $w=(w\wedge x)\vee \tau$,
for some atom $\tau \not\le x$. Since $w^\bperp\cap x=\{\alpha\in
x:\ \omega(\alpha,\beta)=0\quad \forall \ \beta\in w\}$ and for
all $\alpha\in x$, $\omega(\alpha,\beta)=0$ for all $\beta\in
w\wedge x$ (because $x$ is isotropic), then we have that
$$w^\bperp\cap x=\{\alpha\in x :\
\omega(\alpha,\tau)=0\}=\tau^\bperp\cap x .$$ Letting $V$ be the
underlying space of definition \ref{dpg} we have:
\begin{eqnarray*}
dim(w^\bperp\cap x)&=&dim(\tau^\bperp\cap x)\\
&=&dim(\tau^\bperp)+dim(x)-dim(span\{\tau^\bperp,x\})\\
&=&dim V-1+d-dim(span\{\tau^\bperp,x\})\\
\end{eqnarray*}

But $dim(span\{\tau^\bperp,x\})$ can be either
$dim(\tau^\bperp)=dim(V)-1$, or else is $dim(V)$. In the first
case, we would have that $x\subset \tau^\bperp$, which would imply
that $span\{x,\tau\}$ is isotropic, absurd since $x$ is maximal
isotropic and $\tau \not\le x$. Therefore, we have the second
case, and $dim(w^\bperp\cap x)=(dim V)-1+d-dim(V)=d-1$. This
implies that $\theta(x)=(d-1)-(j-1)=d-j$ and $S(x)=\qbin
{d-j}1=R(x)$.

\end{enumerate}
\eop

\begin{lemma}\label{mui}
Consider the same hypothesis of the previous Proposition and let
$\Lambda_j$ be defined as in \ref{deflambdaj}.

Then, there are constants $\mu_0>\mu_1>...>\mu_{d}$ such that for
every $v\in \Lambda_j$, there exists $v'\in \Lambda_{j-1}$ with
$\L(v)=\mu_j v+v'$.
\end{lemma}
\proof It is enough  to show the lemma for each element of the
spanning set $\{\iota(x)\}_{x \in \Omega_j}$. We prove it
inductively starting at $j={d}$.

Taking $\iota(x)\in \Lambda_{d}$, by Lemma \ref{lemasumaiotas}, we
have that $\L(\iota(x))=-\qbin d1\iota(x)+\sum_{z\lcunor
x}\iota(z)$. Since  $\sum_{z\lcunor x}\iota(z) \in
\Lambda_{{d}-1}$ the proposition holds for $j=d$ with
$\mu_{d}=-\qbin d1$.

Assume now that the following inductive hypothesis is true for
$j+1$: $$\mbox{ {\it There are constants}}\ \mu_{j+1}\ \mbox{{\it
such that}}\ \L(\iota(u))=\mu_{j+1}\iota(u)+\sum_{z\lcud
u}\iota(z)\ \forall\ u\in \Omega_{j+1}$$

  Let $\iota(w)\in \Lambda_j$.
  Let's recall that by Lemma \ref{lemawstar} $\iota(w)={1\over \qbin {d-j}1}\sum_{u\gcud w}\iota(u)$.

\vbox{
\begin{eqnarray*}
\L(\iota(w))&=&{1\over \qbin {d-j}1}\sum_{u\gcud w}\L(\iota(u))\\
&=&{1\over \qbin {d-j}1}\sum_{u\gcud w}\Big(\mu_{j+1}\iota(u)+
\sum_{z\lcud u}\iota(z)\Big) \\
&=&\mu_{j+1}{1\over \qbin {d-j}1}\sum_{u\gcud w}\iota(u)
+{1\over \qbin {d-j}1}\sum_{u\gcud w}\sum_{z\lcud u}\iota(z)\\
&=&\mu_{j+1}\iota(w)+
\DDelta\iota(w)+\sum_{v\lcud w}\iota(v)\\
&=&\left(\mu_{j+1}+q^j+q^{d+e-j-1}\right)\iota(w)+\sum_{v\lcud w}\iota(v)\\
\end{eqnarray*}
}

This proves the inductive hypothesis by taking
$\mu_j=\mu_{j+1}+\DDelta$

In particular we have $\mu_j>\mu_{j+1}$ and the lemma is proved.
\eop

\begin{corollary}\label{Linv}
For $j=0,...,d$, $\Lambda_j$ are $\L$-invariant subspaces of
$\RR^{X}$.
\end{corollary}
\proof
 This follows directly by the previous lemma and Corollary \ref{corolambdas}.
\eop

\begin{theorem}\noindent

For $j =0,...,d$, $V_j=\Lambda_j \cap \Lambda_{j-1}^{\perp}$, are
eigenspaces of $\L$ with corresponding eigenvalue $\mu_j= q^{e}
\qbin {d-j}1-\qbin {j}1$.
\end{theorem}
\proof

Take $v\in V_j\ (\subseteq \Lambda_j)$. By Lemma \ref{mui}
$\L(v)=\mu_jv+v'$ with  $v'\in \Lambda_{j-1}$  and by Corollary
\ref{Linv} $\L(v')\in \Lambda_{j-1}$. Then by definition of $V_j$:
\begin{eqnarray*} 0&=&<v,\L(v')>\\
&=&<\L(v),v'>\\
&=&<\mu_jv+v',v'>\\
&=&<\mu_jv,v'>+<v',v'>\\
&=&||v'||^2
\end{eqnarray*}
 thus  $\L(v)=\mu_jv \ \forall \ v\in V_j$.
 Therefore, $\RR^{X}=\bigoplus_{j=0}^{d} V_j$ where each $V_j$ is either zero or an eigenspace of $\L$.

Since $X=\Omega_{d}$ is the set of vertices of a dual polar graph
$\Gamma=(X,E)$ (a distance regular graph of diameter ${d}$),
recall that there are exactly $d+1$ eigenspaces of the adjacency
matrix $A_1$, therefore of $\L$. Thus each $V_j$ is indeed an
eigenspace of $\L$ (hence $V_j\neq 0 \ \forall \ j$) and $\mu_j$
are the eigenvalues of $\L$. Their values can be obtained by
induction from the formula $\mu_j=\mu_{j+1}+\DDelta$, or directly
from Theorem 9.4.3 of \cite{BCN} (since $\mu_j>\mu_{j+1}$).

\eop

\section{Tight Frames for the eigenspaces of a dual polar graph}

In this section we will consider $\Gamma=(X,E)$ be a dual polar
graph of diameter $d$, $\La$ the associated lattice described in
Section \ref{ldpg} and
$$\RR^{X}=\oplus_{j=0}^{d} V_j$$ the corresponding decomposition.

We will give a finite tight frame on each $V_j$ and a formula for
the corresponding constants. We will compute  explicitely the
constant associated to the eigenspace of the second largest
eigenvalue.

\begin{definition}\noindent

Given a vector space $(V,<;>)$ a finite tight frame on $V$ is a
finite set $F \subseteq V$ which satisfies the following
condition: there exists a non-zero constant $\lambda$ such that:
$$ \sum_{v\in F} | <f,v>|^{2} = \lambda \ \| \ f \|^{2}  \quad \forall \ f \in V$$
\end{definition}
As a consequence, $f$ can be expanded as follows: $ f=
\frac{1}{\lambda}\sum_{v\in F} <f,v> v$

\begin{definition}
For $j=0,1,...,d$, let  $U^{j} \in \RR^{X\times X}$ be the matrix
$$(U^{j})_{x,y}=(x,y)^j\quad \mbox{where}\quad (x,y)^j=\sum_{u\in \Omega_j}\iota(u)(x)\iota(u)(y)$$
\end{definition}

\begin{lemma}\label{Uj}
For $j=0,1,...,d$ $$U^{j}=\sum_{l=j}^{d}\qbin l j A_{d-l}$$ where
$A_{i}$ is the $i$-th adjacency matrix of the dual polar graph.

\end{lemma}

 \proof

Let $(x,y) \in X\times X$ and $l=\rk(x\wedge y )$
\[\begin{array}{rcll}
U^{j}_{x,y}&=&(x,y)^j \\
&=&\sum_{u\in \Omega_j}\iota(u)(x)\iota(u)(y)\\
&=&\sum_{u\in \Omega_j}[u\le x][u\le y]\\
&=&\sum_{u\in\Omega_j}[u\le x\wedge y]\\

&=&\left\{\begin{array}{ccl}
 |\{u\in\Omega_j:u \le x\wedge y \}| \ &\mbox{if} &\ j\le l \\
0 &\mbox{if}& \ j > l \end{array} \right. \\
\\
&=& \qbin l j A_{d-l}(x,y)
\end{array}\]
\eop

\begin{definition}\noindent

For $j=0,1,...,d$, let  $\pi_j$ be the orthogonal projection
\\ $\pi_j:\RR^{X}\rightarrow V_j$.
Then for each $u\in \Omega_j$, denote
 $\check u=\pi_j(\iota(u)).$
\end{definition}

Using the previous lemma, we obtain the following:
\begin{corollary} (of Lemma \ref{Uj}) \label{Ueigen} \noindent
For every $j=0,...,d$, $\check u$ is an eigenvector of $U^{j} $
with eigenvalue $\lambda_j=\sum_{l=j}^{d}\qbin l j p_{d-l}(j)$,
where $p_{i}(j)$ are the eigenvalues of $A_i$ corresponding to the
eigenspace $V_j$.
\end{corollary}
\proof  Since $\check u \in V_j$ is an eigenvector of the
adjacency matrices, the Corollary follows from the expression of
$U^j$ given in Lemma \ref{Uj}. \eop

\begin{remark} {\rm Computation of $\lambda_j$.}\label{remlambaj}

From previous Corollary and making the change of variable $i=d-l$,
we have that $\lambda_j=\sum_{i=0}^{d-j}\qbin {d-i} j p_{i}(j)$.

According to pages 261-265, 303-304 of (\cite{BI});
\begin{eqnarray*}
p_{i}(j)&=&u_i(\theta_j) k_i  \qquad \mbox{\rm {where}}\  k_i=p_{ii}^{0}\ \mbox{\rm {and}}\  \ u_i(\theta_j) \  \mbox{\rm is given by the basic hypergeometric series}  \\
u_i(\theta_j)&=& _4\varphi_3 \left(
\begin{array}{llllll}
q^{-i} ,& 0,& q^{-j} ,& -q^{-d-e+j} &  \\
 & & & &;q &,q \\
q^{-d},&0,&0,& & & \end{array}\right)
\mbox{\rm {defined by}} \\
&=& \sum_{t=0}^{\infty}  \frac{(q^{-i};q)_t \ (q^{-j};q)_t \ (-q^{-d-e+j};q)_t \ q^{t}}{(q^{-d};q)_t \ (q;q)_t}\quad  \mbox{\rm {where}}\\
(a;q)_t&=& \left\{\begin{array}{lll} (1-a)...(1-aq^{t-1}) \ &
(t=1,2,...)\\
1 & (t=0 )
\end{array}\right.\\
\mbox{From}&(3)&\mbox{of page 1 and Theorem 9.4.3 of page 275 of
\cite{BCN} } \\k_i&=&\qbin d i q^{\frac{(i^{2}-i)}{2}+i e}
\end{eqnarray*}
\textbf{Caution:} The parameter $"e"$ for dual polar graphs in
page 303 of \cite{BI}, has been replaced by $"e-1"$ to follow the
notation of \cite{BCN} used in Lemma \ref{BCN}.
\end{remark}
\begin{proposition}\label{lambda1}
$$\lambda_1=q^{d-1}\prod_{j=-1}^{d-3} (1+q^{j+e})= q^{d-1} (1+q^{e-1}) \ a_2$$
($a_2$ as in  Lemma \ref{al})
\end{proposition}
\proof
 From Remark \ref{remlambaj}
\begin{eqnarray*}
\lambda_1&=& \sum_{i=0}^{d-1}\qbin {d-i} 1 p_{i}(1)\\
 &=& \sum_{i=0}^{d-1} \qbin {d-i} 1  \ \left( \sum_{t=0}^{1}  \frac{(q^{-i};q)_t \ (q^{-1};q)_t \ (-q^{-d-e+1};q)_t \ q^{t}}{(q^{-d};q)_t \ (q;q)_t}\right)\  \qbin d i
q^{\frac{(i^{2}-i)}{2}+i e}\\
 &=& \sum_{i=0}^{d-1}
\left(1+\frac{(1-q^{-i})(1-q^{-1})(1+q^{-d-e+1}) \
q}{(1-q^{-d})(1-q)}\right) \qbin {d-i} 1 \qbin d i
q^{\frac{(i^{2}-i)}{2}+ie}
\end{eqnarray*}

\vspace{.5em}

since$\qbin {i} 1=\frac{q^{i}-1}{q-1}$ and $\qbin {i}
j=\frac{\qbin {i} 1 \qbin {i-1} 1 ... \qbin {i-j+1} 1}{\qbin {j} 1
\qbin {j-1} 1 ... \qbin {1} 1}$  we have $\qbin {d-i} 1 \qbin d i=
\qbin {d} 1 \qbin {d-1} i$, thus

\vspace{.5em}

\begin{eqnarray*}
 &=& \qbin {d} 1 \sum_{i=0}^{d-1}
\left(1-\frac{\qbin {i} 1 (1+q^{d+e-1}) q^{1-e-i}}{\qbin {d}
1}\right) \qbin {d-1} i
q^{\frac{(i^{2}-i)}{2}+ie} \\
&=& \qbin {d} 1 \sum_{i=0}^{d-1} \qbin {d-1} i
q^{\frac{(i^{2}-i)}{2}+ie}- \sum_{i=0}^{d-1} (1+q^{d+e-1})
q^{1-e-i} \qbin {i} 1 \qbin {d-1} i
q^{\frac{(i^{2}-i)}{2}+ie} \\
&=&   \qbin {d} 1 \sum_{i=0}^{d-1}  \qbin {d-1} i
q^{\frac{(i^{2}-i)}{2}} q^{ie} - \qbin {d-1} 1
(1+q^{d+e-1})\sum_{i=1}^{d-1} \qbin
{d-2}{i-1}q^{1-i} q^{\frac{(i^{2}-i)}{2}}q^{(i-1)e}\\
\end{eqnarray*}
making change of variables $j=i-1$
\begin{eqnarray*}
\lambda_1&=&  \qbin {d} 1 \sum_{i=0}^{d-1}  \qbin {d-1} i
q^{\frac{(i^{2}-i)}{2}} q^{ie} - \qbin {d-1} 1
(1+q^{d+e-1})\sum_{j=0}^{d-2} \qbin {d-2}{j}
q^{\frac{(j^{2}-j)}{2}}q^{je}
\end{eqnarray*}
Using Newton's formula for Gaussian binomials$$\prod_{k=0}^{n-1}
\left(1+q^{k}t \right)= \sum_{k=0}^{n} q^{\frac{(k^{2}-k)}{2}}
\qbin {n} k \ t^{k}$$ with $t=q^{e}$, we have
\begin{eqnarray*}
&=&  \qbin {d} 1  \prod_{j=0}^{d-2} (1+q^{j+e})- \qbin {d-1} 1
(1+q^{d+e-1}) \prod_{j=0}^{d-3} (1+q^{j+e})\\
&=& \left((q^{d-1}+\qbin {d-1} 1)(1+q^{d-2+e})- \qbin {d-1} 1
(1+q^{d+e-1})
\right) \prod_{j=0}^{d-3} (1+q^{j+e}) \\
&=&q^{d-1} (1+q^{e-1})\prod_{j=0}^{d-3} (1+q^{j+e})\\
&=&q^{d-1} (1+q^{e-1}) \ a_2
\end{eqnarray*}\eop

\begin{theorem}\label{tightframe}
For $j=0,...,d$, the set  $\{\check u\}_{u\in\Omega_j}$ is a
finite tight frame for $V_j$, i.e., there exists $c > 0 $ such
that for all $f\in V_j$
$$\sum_{u\in\Omega_j}<\check u,f>\check u= c\ f$$
Moreover, $c=\lambda_j$ is the eigenvalue of $U^j$ corresponding
to any   $\check u$ with $u\in\Omega_j$.
\end{theorem}

\proof

From Corollary \ref{Ueigen}, we have that $\lambda_j \check
u(x)=\sum_{y\in X}(x,y)^j\check u(y)$. Then, for an arbitrary
$v\in \Omega_j$, we have:
\begin{eqnarray*}
<\lambda_j \check u,\check v>&=&\sum_{x\in X}\lambda_j \check u(x)\check v(x)\\
&=&\sum_{x\in X}\left(\sum_{y\in X}(x,y)^j\check u(y)\right)\check v(x) \\
&=&\sum_{x,y\in X}\sum_{w\in \Omega_j}\iota(w)(x)\iota(w)(y)\check u(y)\check v(x)\\
&=&\sum_{w\in \Omega_j}\left(\sum_{y\in X}\iota(w)(y)\check u(y)\right)\left(\sum_{x\in X}\iota(w)(x)\check v(x)\right) \\
&=&\sum_{w\in \Omega_j}<\iota(w),\check u><\iota(w),\check v>\\
&=&\sum_{w\in \Omega_j}<\check w,\check u><\check w,\check v>\qquad ({\rm by\ orthogonality})\\
&=&<\sum_{w\in \Omega_j}<\check w,\check u>\check w,\check v>
\end{eqnarray*}
Since this is true for an arbitrary $v\in \Omega_j$, we conclude
that
 $$\lambda_j \check u=\sum_{w\in \Omega_j}<\check w,\check u>\check w$$
The statement of the theorem follows from the fact that $\{\check
w\}_{w\in \Omega_j}$ span $V_j$. \eop

\begin{lemma}\noindent \label{taucheck}

For all $\tau\in\Omega_1$, \ $\cht=\iota(\tau)-{a_1\over
| X |}\delta_X $ with $a_1$  given in Lemma \ref{al} and \\
$\delta_X:=\iota(\hat{0})=1 \in \RR^{X}$
\end{lemma}
\proof Recall that $<\iota(\hat{0})>=\Lambda_0 \ \subseteq
\Lambda_1=<\{\iota(\tau)\}_{\tau \in \Omega_1}>$, and
$\Lambda_1=\Lambda_0\oplus V_1$. Since $\cht=\pi_1(\iota(\tau))
\in V_1$, we
have $\cht= \iota(\tau) -t.\delta_X$ for some $t \in \RR$.\\
From the fact that $<\cht,\delta_X>=0$ we conclude $t={<
\iota(\tau) ,\delta_X>\over ||\delta_X ||^2}={\sum_{x\in
X}[\tau\subseteq x]\over |X|}=\nn$. \eop

\begin{corollary}\label{pitheorem}
Let $h \in \RR^{X}$, then
$$\pi_1(h)=\sum_{\tau\in \Omega_1}{< \iota(\tau) ;h>\over \lambda_1}\check \tau$$
\end{corollary}

\proof By Theorem \ref{tightframe}, since $\pi_1(h) \in V_1$, we
have
 $$\pi_1(h)
=\sum_{\tau\in \Omega_1}{< \check \tau ;\pi_1(h)>\over
\lambda_1}\check \tau$$ but since
$h=\pi_0(h)+\pi_1(h)+...+\pi_d(h)$, \ $\pi_i(h)\in V_i $\ and \
$<V_i,V_j>=0 \ \forall \ i\neq j $
\begin{eqnarray*} \pi_1(h) &=&\sum_{\tau\in
\Omega_1}{< \check \tau ;h>\over \lambda_1}\check \tau ,\quad \mbox{by lemma} \ \ref{taucheck} \\
&=&\sum_{\tau\in \Omega_1}{< \iota(\tau)-\frac{a_1}{| X |} \delta
;h>\over \lambda_1}\check \tau\\
&=&\sum_{\tau\in \Omega_1}{< \iota(\tau) ;h>\over \lambda_1}\check
\tau -{<\frac{a_1}{| X |} \delta ;h>\over \lambda_1} \sum_{\tau\in
\Omega_1}\check \tau\\
&&\qquad \mbox{but since}\ \sum_{\tau\in \Omega_1} \tau\in \Lambda_0: \\
&=&\sum_{\tau\in \Omega_1}{< \iota(\tau) ;h>\over \lambda_1}\check
\tau
\end{eqnarray*}

 \eop

\section{Application: Norton product on $V_1$}

Given the decomposition $\RR^{X}=V_0\oplus V_1\oplus...\oplus
V_d$, in this section we describe the product of a Norton algebra
attached to the eigenspace $V_1$.

\begin{definition}
The  Norton algebra on $V_1$ is the algebra given by the product
$f\star g=\pi_1(fg)$ for $f,g\in V_1$.

\end{definition}

We want to compute the $\star $ product in $V_1$. Since
$\Lambda_1=span\{\iota(\tau):\tau\in \Omega_1 \}$ the set
$\{\check \tau\}_{\tau\in \Omega_1}$ spans $V_1$.

We want to be able to compute $\cht\star \chs$ in terms of this
set of generators.

\begin{lemma}\label{taustarsigma}
$$\cht\star\chs=
\pi_1( \iota(\tau \vee \sigma))-\nn (\cht+\chs)
$$
\end{lemma}
\proof Recall that $\delta_X=1 \in \RR^{X}$ and that by Lemma
\ref{taucheck} \ $\cht= \iota(\tau) -{a_1\over | X |}\delta_X$.

Observe that $\delta_X$ is the identity of the product of
functions. Then
\begin{eqnarray*}
\cht \star \chs&=&
( \iota(\tau) -{a_1\over | X |}\delta_X)\star (\iota(\sigma)-{a_1\over | X |}\delta_X )\\
&=&\pi_1\big( ( \iota(\tau) -{a_1\over | X |}\delta_X)(\iota(\sigma)-{a_1\over | X |}\delta_X ) \big)\\
&=&\pi_1\big(  \iota(\tau)  \iota(\sigma)-{a_1\over | X |}(  \iota(\tau)  +\iota(\sigma) )+({a_1\over | X |})^2\delta_X \big)\\
&=&\pi_1\big(  \iota(\tau)  \iota(\sigma) \big)-{a_1\over | X
|}\pi_1(  \iota(\tau)  +\iota(\sigma) )+({a_1\over | X
|})^2\pi_1(\delta_X )\quad (\mbox{by
Lemma}\ \ref{iotahat})\\
&=&\pi_1( \iota(\tau \vee \sigma))-\nn (\cht+\chs)\quad
\end{eqnarray*}

\eop

It is clear that in order to complete the description of the
product $\star $ , we need to be able to calculate
 $\pi_1(\iota(\tau \vee \sigma))=\pi_1(\iota(\tau) \iota(\sigma))$.

By Corollary  \ref{pitheorem}:
\begin{eqnarray*}
\pi_1(\iota(\tau \vee
\sigma))&=&\sum_{\rho\in\Omega_1}{<\iota(\rho);\ (\iota(\tau \vee
\sigma)) >\over \lambda_1 }\chr
\end{eqnarray*}

Therefore we need to compute $<\iota(\rho); \iota(\tau \vee
\sigma) >$. We do this in the following:
\begin{lemma}\label{rhotausigma}
$$<\iota(\rho); \iota(\tau \vee \sigma)
>=a_{\rk(\rho\vee\tau\vee\sigma)}$$ where $a_j$ are as in Lemma
\ref{al}.
\end{lemma}
\proof

\begin{eqnarray*}
<\iota(\rho); \iota(\tau \vee \sigma) >&=&
\sum_{x\in X}\iota(\rho)(x) \iota(\tau \vee \sigma) (x)\\
&=&\sum_{x\in X}[\rho\leq x][\tau \vee \sigma\leq x]\\
&=&\sum_{x\in X}[\rho\vee\tau \vee\sigma\leq x]\\
&=&|\{x\in X:\rho\vee\tau \vee\sigma\leq x\}|\\
&=&a_{\rk(\rho\vee\tau\vee\sigma)}\\
\end{eqnarray*}
\eop
\begin{definition}
Given $\tau,\sigma\in \Omega_1$, let:
$$\Psi_j=\{\rho\in \Omega_1: \rk(\rho\vee\tau\vee\sigma)=j\}$$
\end{definition}

\begin{theorem}
$$\cht\star\chs+\nn (\cht+\chs)=\begin{cases}\cht& {\rm if}\ \tau=\sigma\cr
0& {\rm if} \  \tau\vee \sigma=\hat{1} \cr
(1+q^{d-3+e})\sum_{\rho\in\Psi_2}\chr\ +\ \sum_{\rho\in\Psi_3}\chr
\over q^{d-1} (1+q^{e-1})(1+q^{d-3+e}) & {\rm otherwise}  \cr
\end{cases}$$
\end{theorem}

\proof \noindent

By Lemma \ref{taustarsigma}, $\cht\star\chs+\nn
(\cht+\chs)=\pi_1(\iota(\tau \vee \sigma)) $.  Then, the cases
$\tau=\sigma$ and $\tau\vee \sigma =\hat{1}$ follow.

In the remaining case, $\tau\vee \sigma \in \Omega_2$, so

\begin{eqnarray*}
\pi_1( \iota(\tau \vee \sigma) ) &=&
\sum_{\rho\in\Omega_1}{<\iota(\rho); \iota(\tau \vee \sigma) >\over  \lambda_1}\chr \\
&=&\sum_{\rho\in\Omega_1}{a_{\rk(\rho\vee\tau\vee\sigma)}\over
 \lambda_1} \chr
\quad (\diamond)\\
\end{eqnarray*}

Observe that ${\rm rk}(\rho\vee\tau\vee\sigma) \in \{2,3,d+1\}$.
If ${\rm rk} (\rho\vee\tau\vee\sigma)=d+1, a_{d+1}=0$, otherwise
$\rho\vee\tau\vee\sigma\in \Psi_2$ or $\Psi_3$.

Recall that $ \lambda_1=q^{d-1}(1+q^{e-1}) a_2$ and
$a_2=(1+q^{e+d-3}) a_3$ (Remark \ref{remal} and Proposition
\ref{lambda1}).

So $(\diamond)$ becomes
\begin{eqnarray*}
\pi_1( \iota(\tau \vee \sigma) ) &=& {a_{2}\sum_{\rho\in\Psi_2}
\chr +a_{3}\sum_{\rho\in\Psi_3} \chr \over
 \lambda_1}  \\
&=& {(1+q^{e+d-3})\sum_{\rho\in\Psi_2} \chr +\sum_{\rho\in\Psi_3}
\chr \over q^{d-1} (1+q^{e-1})(1+q^{d-3+e})}
\end{eqnarray*}
\eop

\begin{remark}
Similar results are valid for the Johnson, Grassmann and Hamming
cases. They are part of a current research.

\vspace{2em}

 In the Johnson case $J(n,k )$, when $3\le k< {n\over 2}$, we obtain:
$$\cht\star\chs=\begin{cases}(1-2{k\over n} )\cht& {\rm
if}\ \tau=\sigma \cr {2k-n\over n(n-2)}(\cht+\chs)& {\rm if}\
\tau\ne \sigma\end{cases}$$ In the Hamming case,  our formula for
the Norton product reduces to zero. This is also direct from
Theorem 5.2 of \cite{CGS} since $q_{1,1}^1=0$.

In the Grassmann $J_q(n,k)$ case, when $3\le k< {n\over 2}$, we
obtain:

$$\cht\star\chs+{q^k-1\over q^n-1}(\cht+\chs)=\begin{cases} \cht& {\rm
if}\ \tau=\sigma\cr \cr {q^{k-1}-1\over
q(q^{n-2}-1)}\sum_{\rho\in\Psi_2}\chr& {\rm if}\  \tau\neq \sigma
\end{cases}$$
\end{remark}


\begin{thebibliography}{Carolina}

\bibitem {BCN}Brouwer, A. E.; Cohen, A.; Neumaier, A.
{\sl Distance-regular graphs.} Ergebnisse der Mathematik und ihrer
Grenzgebiete. 3. Folge, 18. Berlin etc.: Springer-Verlag. xvii,
495 p. (1989).
\bibitem {BI}  Bannai E.; Ito T.{\sl Algebraic Combinatorics I:Association Schemes.} Benjamin Cummings. London,1984
\bibitem {BF}Benedetto, John J.; Fickus, Matthew  {\sl Finite normalized tight
frames.} Adv. Comput. Math. \textbf{18}, No.2-4, 357-385 (2003).
\bibitem {CGS} Cameron, P.; Goethals,J.; Seidel, J.,
{\sl The Krein condition, spherical designs, Norton algebras and
permutation groups},  Proc. Kon. Nederl. Akad. Wetensch. (A)
\textbf{81} (1978),  196-206.
\bibitem {DR} Diaconis, Persi; Rockmore, Daniel {\sl Efficient computation
of isotypic projections for the symmetric group.} DIMACS, Ser.
Discrete Math. Theor. Comput. Sci. \textbf{11}, 87-104 (1993).
\bibitem {FMRHO} Foote, Richard; Mirchandani, Gagan; Rockmore, Daniel N.; Healy,
Dennis; Olson, Tim {\sl A wreath product group approach to signal
and image processing. I: Multiresolution analysis.} IEEE Trans.
Signal Process. \textbf{48}, No.1, 102-132 (2000).
\bibitem {KC1} Kova\v{c}evi\'c, J.; Chebira, A., {\sl Life beyond
bases: The advent of frames (Part I)}, IEEE Signal Processing
Mag., \textbf{24}, no. 4 (2007), pp. 86-104.
\bibitem {KC2}  Kova\v{c}evi\'c, J.; Chebira,  A., {\sl Life beyond
bases: The advent of frames (Part II)}, IEEE SP Mag., \textbf{24},
no. 5 (2007), pp. 115-125.
\bibitem {MFRHO} Mirchandani, Gagan; Foote, Richard; Rockmore, Daniel N.; Healy,
Dennis; Olson, Tim {\sl A wreath product group approach to signal
and image processing. II: Convolution, correlation, and
applications.} IEEE Trans. Signal Process. \textbf{48}, No.3,
749-767 (2000).
\bibitem{S1} Stanton, Dennis. {\sl Some q-Krawtchouk polynomials on Chevalley
groups.} Am. J. Math. \textbf{102}, 625-662 (1980).
\bibitem{S2} Stanton, Dennis. {\sl Orthogonal polynomials and Chevalley groups.}
Special functions: Group theoretical aspects and applications,
Math. Appl., D. Reidel Publ. Co. \textbf{18}, 87-128 (1984).
\bibitem{VW1}Vale, Richard; Waldron, Shayne.  {\sl Tight frames and their
symmetries.}
 Constructive Approximation \textbf{21}, No. 1, 83-112 (2005).
\bibitem{VW2} Vale, Richard; Waldron, Shayne.  {\sl Tight frames generated by finite
nonabelian groups.}  Numer. Algorithms \textbf{48}, No. 1-3, 11-27
(2008).
\end{thebibliography}
\end{document}